\documentclass[12pt]{article}
\usepackage{latexsym, amssymb, amsmath, amscd, amsfonts, epsfig, graphicx, colordvi}

\newtheorem{thm}{Theorem}[section]

\def\pf{\noindent{\it Proof.} }
\def\qed{\nopagebreak\hfill{\rule{4pt}{7pt}}
\medbreak}

\def\qed{\nopagebreak\hfill{\rule{4pt}{7pt}}
\medbreak}

\newenvironment{kst}
{\setlength{\leftmargini}{2.4\parindent}
\begin{itemize}
\setlength{\itemsep}{-0.5mm}} {\end{itemize}}

\begin{document}

\begin{center}

{\Large \bf Weighted Forms of Euler's Theorem}\\[15pt]
\end{center}

\begin{center}
{ William Y. C. Chen}$^{1}$ \quad  and  \quad {Kathy Q. Ji}$^{2}$

   Center for Combinatorics, LPMC\\
   Nankai University, Tianjin 300071, P.R. China

   \vskip 1mm

   Email: $^1$chen@nankai.edu.cn, $^2$ji@nankai.edu.cn
\end{center}

\vskip 6mm \noindent {\bf Abstract.}  In answer to a question of
Andrews about finding combinatorial proofs of two identities in
Ramanujan's ``Lost" Notebook, we obtain weighted forms of Euler's
theorem on partitions with odd parts and distinct parts. This work
is inspired by the insight of Andrews on the connection between
Ramanujan's identities and Euler's theorem. Our combinatorial
formulations of Ramanujan's identities rely on the notion of
rooted partitions. Iterated Dyson's map and Sylvester's bijection
are the main ingredients in the weighted forms of Euler's theorem.



\section{ Introduction}

This paper is concerned with the combinatorial treatments of the
   following two identities from Ramanujan's
``Lost" Notebook:
\begin{equation}
\sum_{\tiny n=0}^{\infty}\left[ (-q;q)_{\infty}-(-q;q)_n \right]
=(-q;q)_{\infty} \left[ -\frac{1}{2}+\sum_{d=1}^{\infty}
\frac{q^d}{1-q^d} \right] +\frac{1}{2}\left[1+\sum_{n=1}^{\infty}
\frac{q^{n+1\choose 2}}{(-q;q)_n}\right],\label{eqn1}
\end{equation}
\begin{equation}
\sum_{n=0}^{\infty}\left[\frac{1}{(q;q^2)_{\infty}}-\frac{1}{(q;q^2)_n}\right]
=(-q;q)_{\infty}\left[-\frac{1}{2}+\sum_{d=1}^{\infty}\frac{q^{2d}}{1-q^{2d}}\right]
+\frac{1}{2}\left[1+\sum_{n=1}^{\infty}\frac{q^{n+1\choose
2}}{(-q;q)_n}\right],\label{eqn2}
\end{equation}
where the $q$-shifted factorial is defined by $(x; q)_0 = 1$ and
for $n\geq 1$,
\[
(x; q)_n = (1 - x)(1 - qx) \cdots (1- q^{n-1}x).\]

Andrews \cite{and86} has obtained algebraic proofs of the above
identities by differentiation.
 Furthermore he asked ``Can a `near
bijection' be provided between the weighted counts of partitions
given by the left sides of \eqref{eqn1} and \eqref{eqn2} and the
convolution of partition functions generated by the first
summation of the right sides of \eqref{eqn1} and  \eqref{eqn2}?"
Andrews also gave an insightful remark that these two identities
may be seen as closely related to Euler's result although not
strictly generalizations of it, and pointed out the combinatorial
possibilities of studying weighted counts of partitions such as
related to two identities. Our work is indeed inspired by the idea
of Andrews.

Recently, Andrews, Jim\'{e}nez-Urroz and Ono proved many
identities related to the Dedekind eta-function in \cite{and01},
including the above two identities.   Chapman  \cite{cha02} found
a combinatorial formulation \eqref{eqn1}. But he did not get a
combinatorial correspondence and remarked that it would be
interesting to find one. In this paper, we first obtain
combinatorial formulations of the Ramanujan's identities
(\ref{eqn1}) and (\ref{eqn2}) based on a new interpretation of the
first summation on the right hand side of (\ref{eqn1}). We further
obtain weighted counting theorems for partitions into odd parts
and distinct parts, which can be regarded as weighted forms of
Euler's theorem. Then we establish the connections between the
Ramanujan's identities and our weighted forms of Euler's theorem,
just as anticipated by Andrews \cite{and86}. The weighted forms of
Euler's theorem can be derived combinatorially by using
Sylvester's bijection and iterated Dyson's map.

This paper is organized as follows. We give a brief review of
Sylvester's bijection and iterated Dyson's map in Section 2, and
obtain  weighted forms of Euler's theorem. In Section 3, we
introduce the notion of rooted partitions and obtain generating
functions for rooted partitions as well as identities on rooted
partitions. In Section 4, we establish the connections between
weighted forms of Euler's theorem (Theorems \ref{thm1} and
\ref{thm2}) and Ramanujan's identities (\ref{eqn1}) and
\eqref{eqn2} via identities on rooted partitions.

\section{ Sylvester's Bijection and Iterated Dyson's Map}

In this section, we give  several weighted forms of Euler's
theorem from Sylvester's bijection and iterated Dyson's map. We
first recall some terminology on partitions  in \cite{And76}. {\it
A partition} $\lambda$ of a positive integer $n$ is a finite
nonincreasing sequence of positive integers
$\lambda_1,\,\lambda_2,\ldots,\,\lambda_r$ such that
$\sum_{i=1}^r\lambda_i=n.$  Then $\lambda_i$ are called the parts
of $\lambda$, $\lambda_1$ is its largest part. The number of parts
of  $\lambda$ is called the length of $\lambda$ denoted by
$l(\lambda).$ Let $n_{\lambda}(d)$ be the number of parts equal to
$d$ in $\lambda$, then we have $l(\lambda)=\sum_d n_{\lambda}(d).$
The weight of $\lambda$ is the sum of  parts of $\lambda$, denoted
by $|\lambda|.$

{\it The rank of a partition} $\lambda$ introduced  by Dyson
\cite{dys44} is defined as the largest part minus the number of
parts, which is usually denoted by $r(\lambda).$  As a convention,
we shall assume that the empty partition  has  rank zero. For a
partition $\lambda=(\lambda_1,\ldots,\,\lambda_r),$ we define the
conjugate partition
$\lambda'=(\lambda_1',\lambda_2',\ldots,\,\lambda'_t)$ of
$\lambda$ by setting $\lambda'_i$ to be the number of parts of
$\lambda$ that are greater than or equal to $i$. Clearly, we have
$l(\lambda)=\lambda'_1$ and $\lambda_1=l(\lambda').$

The set of  partitions of $n$ into distinct parts is denoted by
$D_n$, and the set of  partitions of $n$ into odd parts  is
denoted by $O_n.$ Euler's theorem states that $|D_n|=|O_n|$ for
$n\geq 1$, which follows from the following generating function
identity:
$$(-q;q)_{\infty}=\frac{1}{(q;q^2)_{\infty}}.$$
  Sylvester's bijection \cite{syl82} and iterated Dyson's map\cite{And83} are two
correspondences between $D_n$ and $O_n$. As we will see, they play
a key role in the proofs of the weighted forms of Euler's theorem.

 There are several ways to describe Sylvester bijection
\cite{And84,bes94, bre99,  mac84}. Here we give a description by
using 2-modular diagram as given by Bessenrodt \cite{bes94}.

{\bf Sylvester's bijection $\varphi$}:  Given a partition
$\lambda$ of $n$ with odd parts, represent each part $2m+1$ by a
row of $m$ 2's and a 1 at the end. This diagram is called the
2-modular diagram of $\lambda$.  Decompose the 2-modular diagram
into hooks $H_1,\,H_2,\ldots$ with the diagonal boxes as corners.
Let $\mu_1$ be the number of squares in $H_1$, let $\mu_2$ be the
number of 2's in $H_1$, let $\mu_3$ be the number of squares in
$H_2$, let $\mu_4$ be the number of 2's in $H_2$, and so on. Set
$\varphi(\lambda)=\mu=(\mu_1,\, \mu_2,\,\mu_3,\, \ldots)$. Then
$\varphi(\lambda)$ is clearly a partition with distinct parts, see
Figure $1$.

{\bf The inverse map $\varphi^{-1}$}: Let
$\mu=(\mu_1,\,\mu_2,\ldots,\,\mu_{2k-1},\,\mu_{2k})$ be a
partition of $n$ into distinct parts, where $\mu_i >0$ for $1 \leq
i \leq 2k-1$ and $\mu_{2k} \geq 0.$ First we consider the part
$\mu_{2k}$, and write down $\mu_{2k}$ 2's in a row and add a 1 to
 the end of the first row, then add $(\mu_{2k-1}-\mu_{2k}-1)$
1's to the first column. Let us denote this hook by $H_k$. Note
that the 2's can only appear in the first row in this hook. Let us
continue to consider the parts $\mu_{2k-3},\, \mu_{2k-2}$. The
hook $H_{k-1}$ is constructed as follows. There will be
$\mu_{2k-2}$ 2's and $\mu_{2k-3}-\mu_{2k-2}$ 1's in $H_{k-1}$. If
there is a 1 in the $i$-th of $H_k$, then there must a 2 on the
left of the 1 in $H_k$. The rest of the 2's will have to be put in
the first row of $H_{k-1}$. Then the 1's are easily dispatched in
the first row and the first column.  Now we may repeat the above
procedure to construct a partition with odd parts.

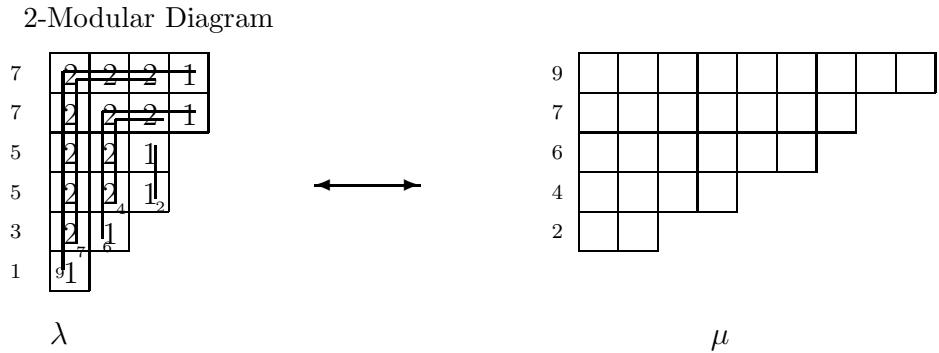
\begin{figure}[h]
\begin{center}
\begin{picture}(100,120)
\put(-100,100){\line(1,0){60}} \put(-100,100){\line(0,-1){90}}
\put(-100,10){\line(1,0){15}} \put(-85,10){\line(0,1){90}}
\put(-100,25){\line(1,0){30}}\put(-70,25){\line(0,1){75}}
\put(-100,40){\line(1,0){45}}\put(-55,40){\line(0,1){60}}
\put(-100,55){\line(1,0){45}}\put(-100,70){\line(1,0){60}}
\put(-40,70){\line(0,1){30}}\put(-100,85){\line(1,0){60}}
\put(-110,110){{\small 2-Modular Diagram}}
\put(-95,88){2}\put(-80,88){2}\put(-65,88){2}\put(-50,88){1}
\put(-95,73){2}\put(-80,73){2}\put(-65,73){2}\put(-50,73){1}
\put(-95,58){2}\put(-80,58){2}\put(-65,58){1} \put(-95,43){2}
\put(-80,43){2}\put(-65,43){1}\put(-95,28){2}\put(-80,28){1}
\put(-95,13){1}\put(-115,90){{\scriptsize7}}
\put(-115,75){{\scriptsize7}}\put(-115,60){{\scriptsize5}}
\put(-115,45){{\scriptsize5}}\put(-115,30){{\scriptsize3}}
\put(-115,15){{\scriptsize1}} {\thicklines
\put(-95,93){\line(1,0){50}}\put(-95,93){\line(0,-1){75}}
\put(-98,15){{\tiny 9}}\put(-90,90){\line(1,0){30}}
\put(-90,90){\line(0,-1){62}}\put(-90,23){{\tiny 7}}
\put(-80,78){\line(1,0){35}}\put(-80,78){\line(0,-1){48}}
\put(-80,25){{\tiny6}}\put(-75,75){\line(1,0){18}}
\put(-75,75){\line(0,-1){32}}\put(-75,39){{\tiny 4}}
\put(-60,65){\line(0,-1){20}}\put(-60,40){{\tiny 2}}
\put(10,50){\vector(1,0){30}}\put(30,50){\vector(-1,0){30}}}
\put(100,100){\line(1,0){135}}\put(100,100){\line(0,-1){75}}
\put(100,25){\line(1,0){30}}\put(130,25){\line(0,1){75}}
\put(100,40){\line(1,0){60}}\put(160,40){\line(0,1){60}}
\put(100,55){\line(1,0){90}}\put(190,55){\line(0,1){45}}
\put(100,70){\line(1,0){105}}\put(205,70){\line(0,1){30}}
\put(100,85){\line(1,0){135}}\put(235,85){\line(0,1){15}}
\put(115,100){\line(0,-1){75}}\put(145,100){\line(0,-1){60}}
\put(175,100){\line(0,-1){45}}\put(220,100){\line(0,-1){15}}
\put(90,90){{\scriptsize9}}\put(90,75){{\scriptsize7}}
\put(90,60){{\scriptsize6}}\put(90,45){{\scriptsize4}}
\put(90,30){{\scriptsize2}} \put(150,-10){$\mu$}
\put(-100,-10){$\lambda$}
\end{picture}
\end{center}
\caption{Sylvester's bijection $\varphi:\ (7,\,7,\,5,\,5,\,3,\,1)
\mapsto (9,\,7,\,6,\,4,\,2).$}
\end{figure}

We now give a brief description of the  bijection due to  Andrews
\cite{And83}, which is called iterated Dyson's map. This
correspondence gives a combinatorial interpretation of a partition
theorem of Fine \cite{fin48, fin88}. Our presentation follows the
survey of Pak \cite{pak03}.

We first describe Dyson's map \cite{dys69}. Denote by $H_{n,\,r}$
and $G_{n,\,r}$ the sets of partitions of $n$ with rank at most
$r$ and at least $r$, respectively. Dyson's map $\psi_r$ is a
bijection between $H_{n,\,r+1}$ and $G_{n+r,\,r-1}$.

{\bf Dyson's map $\psi_r$}: Start with a Young diagram
corresponding to a partition $\lambda \in H_{n,\,r+1}$. Note that
$\lambda$ has $l=l(\lambda)$ parts, where $l(\lambda)$ is the
length, or the number of parts of $\lambda$. Remove the first
column, add  $l+r$ squares to the top row to obtain a Young
diagram, it follows that the resulting Young diagram is a
partition $\mu \in G_{n+r,\,r-1}$. It is easy to see that the
above procedure is reversible. Hence, Dyson's map $\psi_r$ is a
bijection. An example is illustrated in Figure 2.
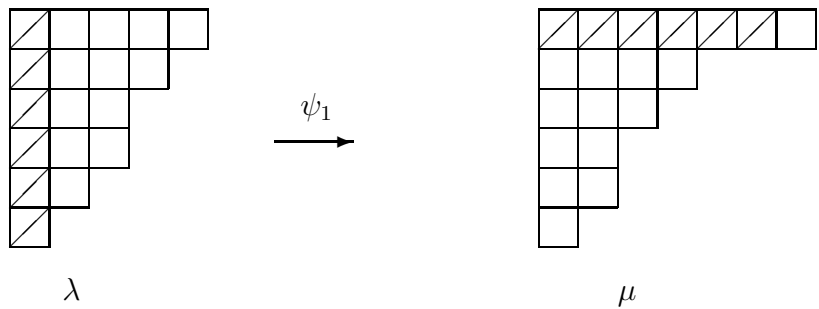
\begin{figure}[h] \label{f2}
\begin{center}
\begin{picture}(100,120)
\put(-100,100){\line(1,0){75}} \put(-100,100){\line(0,-1){90}}
\put(-100,10){\line(1,0){15}} \put(-85,10){\line(0,1){90}}
\put(-100,25){\line(1,0){30}}\put(-70,25){\line(0,1){75}}
\put(-100,40){\line(1,0){45}}\put(-55,40){\line(0,1){60}}
\put(-100,55){\line(1,0){45}}\put(-100,70){\line(1,0){60}}
\put(-40,70){\line(0,1){30}}\put(-100,85){\line(1,0){75}}
\put(-25,100){\line(0,-1){15}}\put(-100,85){\line(1,1){15}}
\put(-100,70){\line(1,1){15}}\put(-100,55){\line(1,1){15}}
\put(-100,40){\line(1,1){15}}\put(-100,25){\line(1,1){15}}
\put(-100,10){\line(1,1){15}}
{\thicklines\put(0,50){\vector(1,0){30}} \put(10,60){$\psi_1$}}
\put(100,100){\line(1,0){105}}\put(100,100){\line(0,-1){90}}
\put(100,10){\line(1,0){15}}\put(115,10){\line(0,1){90}}
\put(100,25){\line(1,0){30}}\put(130,25){\line(0,1){75}}
\put(100,40){\line(1,0){30}}\put(100,55){\line(1,0){45}}
\put(145,55){\line(0,1){45}}\put(100,70){\line(1,0){60}}
\put(160,70){\line(0,1){30}}\put(100,85){\line(1,0){105}}
\put(205,85){\line(0,1){15}}\put(175,100){\line(0,-1){15}}
\put(190,100){\line(0,-1){15}}\put(100,85){\line(1,1){15}}
\put(115,85){\line(1,1){15}}\put(130,85){\line(1,1){15}}
\put(145,85){\line(1,1){15}}\put(160,85){\line(1,1){15}}
\put(175,85){\line(1,1){15}} \put(130,-10){$\mu$}
\put(-80,-10){$\lambda$}

\end{picture}
\end{center}
\caption{$\lambda=(5,\,4,\,3,\,3,\,2,\,1)$ and
$\mu=(7,\,4,\,3,\,2,\,2,\,1).$ }
\end{figure}

We are now ready to describe  iterated Dyson's map $\phi$: $O_n$
$\mapsto$ $D_n.$

{\bf Iterated Dyson's map $\phi$}: Let $\lambda=(\lambda_1,
\lambda_2, \ldots, \lambda_l)$ be a partition of $n$ into odd
parts. We construct a partition $\mu$ of $n$ from $\lambda$ by the
following process. Let $\nu\,^l=(\lambda_l)$ and let $\nu\,^i$
denote the partition obtained by applying Dyson's map
$\psi_{\lambda_i}$ to $\nu\,^{i+1}$ , i.e.
$\nu\,^i=\psi_{\lambda_i}(\nu\,^{i+1})$.  Finally, set
$\mu=\nu\,^1$. Since
$\nu\,^i=\lambda_i+\lambda_{i+1}+\cdots+\lambda_l,$ one sees that
 $|\mu|=|\lambda|$. Furthermore $\mu$ is a partition into distinct parts and
  the iterated Dyson's map $\phi$ is a bijection
 \cite{pak03}.

 The inverse of  iterated Dyson's maps is described as a recursive
procedure.  Let $\mu=(\mu_1,\,\mu_2,\ldots,\,\mu_l)$ be a
partition of $n$ into distinct parts. Set
$\lambda_1=r(\mu)=\mu_1-l(\mu)$ if $r(\mu)$ is odd; otherwise set
$\lambda_1=r(\mu)+1=\mu_1-l(\mu)+1$. Applying the inverse of
Dyson's $\psi_{\lambda_1}^{-1}$ to $\mu,$ we get a partition
$\nu\,^2=\psi_{\lambda_1}^{-1}(\mu).$  Iterating the above
procedure to $\nu\,^j$ $(j=2,\,3,\,4,\ldots)$, we obtain a
partition $\lambda=(\lambda_1,\,\lambda_2,\ldots)$ with odd parts.
 Figure 3 is an illustration of  iterated Dyson's
 map.

\begin{figure}[h]\label{1}
\begin{center}
\setlength{\unitlength}{0.25mm}
\begin{picture}(300,50)
\put(-100,0){\line(1,0){15}}\put(-100,0){\line(0,1){15}}
\put(-100,15){\line(1,0){15}}\put(-85,15){\line(0,-1){15}}
\put(-80,8){\vector(1,0){15}} \put(-55,0){\line(1,0){60}}
\put(5,0){\line(0,1){15}}\put(-55,0){\line(0,1){15}}
\put(-55,15){\line(1,0){60}}\put(-40,15){\line(0,-1){15}}
\put(-25,15){\line(0,-1){15}}\put(-10,15){\line(0,-1){15}}
\put(15,8){\vector(3,1){15}} \put(40,-5){\line(1,0){45}}
\put(40,-5){\line(0,1){30}}\put(40,25){\line(1,0){60}}
\put(100,25){\line(0,-1){15}}\put(100,10){\line(-1,0){60}}
\put(85,-5){\line(0,1){30}}\put(70,-5){\line(0,1){30}}
\put(55,-5){\line(0,1){30}}\put(110,8){\vector(1,0){20}}
\put(140,-10){\line(1,0){30}}\put(170,-10){\line(0,1){45}}
\put(140,-10){\line(0,1){45}}\put(140,35){\line(1,0){105}}
\put(245,35){\line(0,-1){15}}\put(245,20){\line(-1,0){105}}
\put(230,35){\line(0,-1){15}}\put(215,35){\line(0,-1){15}}
\put(200,35){\line(0,-1){15}}\put(185,35){\line(0,-1){15}}
\put(155,35){\line(0,-1){15}}\put(140,5){\line(1,0){45}}
\put(185,5){\line(0,1){30}}\put(155,-10){\line(0,1){30}}
\put(260,8){\vector(1,0){20}}\put(290,-20){\line(0,1){60}}
\put(290,40){\line(1,0){120}}\put(410,40){\line(0,-1){15}}
\put(410,25){\line(-1,0){120}}\put(290,10){\line(1,0){90}}
\put(380,10){\line(0,1){30}}\put(290,-5){\line(1,0){30}}
\put(320,-5){\line(0,1){45}}\put(290,-20){\line(1,0){15}}
\put(305,-20){\line(0,1){60}}\put(335,40){\line(0,-1){30}}
\put(350,40){\line(0,-1){30}}\put(365,40){\line(0,-1){30}}
\put(395,40){\line(0,-1){15}}

\end{picture}
\end{center}
\caption{$\lambda=(5,\,5,\,3,\,3,\,1)$ and $\mu=(8,\,6,\,2,\,1).$}
\end{figure}
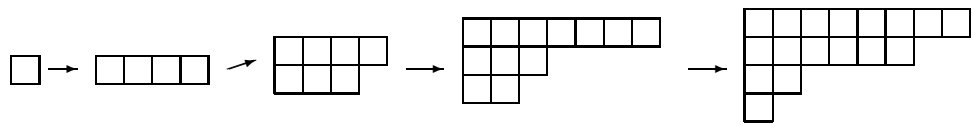

 When applying Sylvester's bijection,
 we see that each partition $\mu$ of $n$ into distinct
 parts with maximal part $\mu_1$ corresponds to a partition $\lambda$
of $n$ into odd parts with the maximal part $\lambda_1$ and the
length $l(\lambda)$ such that
$\mu_1=\frac{\lambda_1-1}{2}+l(\lambda)$ or
$2\mu_1+1=\lambda_1+2l(\lambda).$ Thus we have the following
weighted form of Euler's theorem:

\begin{thm}\label{t1}
The sum of $\mu_1$ {\rm(}{or $2\mu_1+1$\rm)} over all the
partitions $\mu$ of $n$ into distinct parts   equals to the sum of
$\frac{\lambda_1-1}{2}+l(\lambda)$ {\rm(}{or
$\lambda_1+2l(\lambda)$\rm)} over all the partitions $\lambda$  of
$n$ into odd parts, namely,
\begin{equation}\label{euler1}
\sum_{\mu \in D_n}\mu_1=\sum_{\lambda \in
O_n}\left(\frac{\lambda_1-1}{2}+l(\lambda)\right),
\end{equation}
or equivalently,
\begin{equation}\label{euler2}
\sum_{\mu \in D_n}(2\mu_1+1)=\sum_{\lambda \in
O_n}(\lambda_1+2l(\lambda)).
\end{equation}
\end{thm}

 From iterated Dyson's
map, we see that  a  partition $\lambda$ of $n$ into odd parts
with maximal part $\lambda_1$ corresponds to a partition $\mu$  of
$n$ into distinct parts  with  rank $r(\mu)$ such that
\[ r(\mu)+\frac{1+(-1)^{r(\mu)}}{2}=\lambda_1.\]
Thus we obtain the following weighted form of Euler's theorem:

\begin{thm}\label{t2}
The sum of $\mu_1-l(\mu)+\frac{1+(-1)^{r(\mu)}}{2}$ over all
partitions $\mu$ of $n$ into distinct parts   equals  the sum of
$\lambda_1$ over all  partitions $\lambda$ of $n$ into odd parts,
namely,
\begin{equation}\label{euler3}
\sum_{\mu \in
D_n}\left(\mu_1-l(\mu)+\frac{1+(-1)^{r(\mu)}}{2}\right)=\sum_{\lambda
\in O_n}\lambda_1.
\end{equation}
\end{thm}

Now we consider the set of the partitions $\mu$ of $n$ into
distinct parts with multiplicities
$l(\mu)+\mu_1+\frac{1-(-1)^{r(\mu)}}{2}$.  The number of such
partitions of $n$ with the multiplicities taken into account
equals  the number of the elements in the set of partitions of $n$
into distinct parts with multiplicities $2\mu_1+1$ minus the
number of the elements in the set of partitions of $n$ into
distinct parts with multiplicities
$\mu_1-l(\mu)+\frac{1+(-1)^{r(\mu)}}{2}$. According to Theorems
\ref{t1}, \ref{t2}, we obtain the following weighted form of
Euler's theorem which will be used in the  combinatorial proof of
Ramanujan's identity \eqref{eqn1}.

\begin{thm}\label{thm1}
The sum of $l(\mu)+\mu_1+\frac{1-(-1)^{r(\mu)}}{2}$ over all the
partitions $\mu$ of $n$ into distinct parts equal to the sum of
$2l(\lambda)$ over all the partitions $\lambda$ of $n$ into odd
parts, namely,

\begin{equation}\label{we1}
\sum_{\mu \in
D_n}\left(l(\mu)+\mu_1+\frac{1-(-1)^{r(\mu)}}{2}\right)
=\sum_{\lambda \in O_n}2l(\lambda).
\end{equation}

\end{thm}

Next we consider the set of partitions $\mu$ of $n$ into distinct
parts with multiplicities $ l(\mu)+\frac{1-(-1)^{r(\mu)}}{2}$.
 The number of
such partitions with multiplicities equals  the number of elements
in the set of partitions of $n$ into distinct parts with
multiplicities $\mu_1+1$ minus the number of elements in the set
of partitions of $n$ into distinct parts with multiplicities
$\mu_1-l(\mu)+\frac{1+(-1)^{r(\mu)}}{2}$, according to Theorems
\ref{t1}, \ref{t2} and Euler's theorem, we  obtain the following
weighted form of Euler's theorem which will be used in the
combinatorial proof of Ramanujan's identity \eqref{eqn2}.

\begin{thm}\label{thm2}
The sum of $l(\mu)+\frac{1-(-1)^{r(\mu)}}{2}$ over all the
partitions $\mu$  of $n$ into distinct parts    equal to the sum
of $l(\lambda)-\frac{\lambda_1-1}{2}$ over all the partitions
$\lambda$  of $n$ into odd parts, namely,
\begin{equation}\label{we2}
\sum_{\mu \in D_n}\left(l(\mu)+\frac{1-(-1)^{r(\mu)}}{2}\right)
=\sum_{\lambda \in
O_n}\left(l(\lambda)-\frac{\lambda_1-1}{2}\right).
\end{equation}
\end{thm}

\section{ Rooted Partitions}

Inspired by the suggestion of  Andrews  \cite{and86}, we are led
to consider weighted counting of partitions in order to give
combinatorial proofs of Ramanujan's identities \eqref{eqn1} and
\eqref{eqn2}. To this end, we introduce the notion of rooted
partitions which can be regarded as a weighted version of ordinary
partitions. In some sense, rooted partitions are related to
``overpartitions'' (see Corteel and Lovejoy \cite{cor04}) and
``partitions with designated summand'' of Andrews-Lewis-Lovejoy
\cite{and02}.

{\it A rooted partition} of $n$ can be formally defined as a pair
of partitions $(\lambda, \mu)$, where $|\lambda|+|\mu| = n$ and
$\mu$ is a nonempty partition with equal parts.  Intuitively, a
rooted partition is a partition in which some equal parts are
represented as barred elements. The union of the parts of
$\lambda$ and $\mu$ are regarded as the parts of the rooted
partition $(\lambda, \mu)$.

\noindent For example, there are twelve rooted partitions of $4$:
\[
\begin{array}{ll}
 \bar{4},\,\bar{3}+1,\,3+\bar{1},\,
\bar{2}+2,\,\bar{2}+\bar{2},\, \bar{2}+1+1,\,
2+\bar{1}+1,\,2+\bar{1}+\bar{1},\,\bar{1}+1+1+1,& \\[8pt]
\bar{1}+\bar{1}+1+1,\,\bar{1}+\bar{1}+\bar{1}+1,\,
\bar{1}+\bar{1}+\bar{1}+\bar{1}.&
$$
\end{array}
\]

 \noindent There are three rooted partitions of $4$ with
distinct parts: $\bar{4},\,\bar{3}+1,\,3+\bar{1}.$

\noindent There are six rooted partitions of 4 with odd parts:
$$\bar{3}+1,\,3+\bar{1},\,\bar{1}+1+1+1,\,
\bar{1}+\bar{1}+1+1,\,\bar{1}+\bar{1}+\bar{1}+1,
\bar{1}+\bar{1}+\bar{1}+\bar{1}.$$

\noindent A rooted partition $(\lambda, \mu)$ is said to be a
rooted partition with almost distinct parts if $\lambda$ has
distinct parts. There are nine  rooted partitions of $4$ with
almost distinct parts:
$$\bar{4},\,\bar{3}+1,\,3+\bar{1},\,\bar{2}+2,\,\bar{2}+\bar{2},
\,2+\bar{1}+1,\,2+\bar{1}+\bar{1},\,1+\bar{1}+\bar{1}+\bar{1}
,\,\bar{1}+\bar{1}+\bar{1}+\bar{1}.$$

 It is easy to see that the
generating function for rooted partitions with distinct parts
equals
\begin{equation}\label{eqn4}
\sum_{d=1}^{\infty}q^d\prod_{j\neq d}^{\infty}(1+q^j).
\end{equation}
On the other hand, the generating function for rooted partitions
with odd parts equals
\begin{equation}\label{eqn5}
\frac{1}{(q;q^2)_{\infty}}\sum_{d=0}^{\infty}
\frac{q^{2d+1}}{1-q^{2d+1}}.
\end{equation}
The generating function for rooted partitions with almost distinct
parts equals
\begin{equation} \label{g-x}
(-q;q)_{\infty}\sum_{d=1}^{\infty} \frac{q^d}{1-q^d}.
\end{equation}

We now define the {\it root size of a rooted partition}
$(\lambda,\, \mu)$ as the number of parts of $\mu$. Then the
generating function for rooted partitions into almost distinct
parts with even root size equals
\begin{equation} \label{g-y}
(-q;q)_{\infty}\sum_{d=1}^{\infty} \frac{q^{2d}}{1-q^{2d}}.
\end{equation}

We have the following identity on rooted partitions:

\begin{thm}\label{lem2}
The number of the rooted partitions of $n$ into almost distinct
parts with even root size  plus the number of the rooted
partitions of $n$ with distinct parts equals the number of rooted
partitions of $n$ with odd parts.
\end{thm}

We first give a generating function proof of the above theorem.

 \pf
The sum of the two numbers have the following generating function
\begin{align*}
&(-q;q)_{\infty}\sum_{d=1}^{\infty}\frac{q^{2d}}{1-q^{2d} }
+\sum_{d=1}^{\infty}q^d\prod_{n\neq d}^{\infty}(1+q^{n})\\[5pt]
&=(-q;q)_{\infty}\left( \sum_{d=1}^{\infty}\frac{q^{2d}}{1-q^{2d}}
+\sum_{d=1}^{\infty}\frac{q^{d}-q^{2d}}{(1-q^d)(1+q^d)} \right)\\[5pt]
&=(-q;q)_{\infty}\sum_{d=1}^{\infty}\frac{q^d+q^{2d}-q^{2d}}{1-q^{2d}}\\[5pt]
&=(-q;q)_{\infty}\left(\sum_{d=1}^{\infty}\frac{q^d}{1-q^d}
-\sum_{d=1}^{\infty}\frac{q^{2d}}{1-q^{2d}}\right)\\[5pt]
&=\frac{1}{(q;q^2)_{\infty}}
\sum_{d=0}^{\infty}\frac{q^{2d+1}}{1-q^{2d+1}}.
\end{align*}
This implies the desired statement for rooted partitions. \qed

We now present  a combinatorial proof of the above theorem in
terms of an involution and a bijection. We need the following
fact:

\begin{thm} \label{o-1}
 The number of  rooted partitions  of $n$ into almost
distinct parts with odd root size equals  the number of the rooted
partitions of $n$ into almost distinct parts with even root size
plus the number of the rooted partitions of $n$ with distinct
parts.
\end{thm}

\pf We now construct an involution $\tau$ on the set of rooted
partitions of $n$ with almost distinct parts except those strictly
with distinct parts. More precisely, the involution $\tau$ is on
the set of rooted partitions $(\lambda,\, \mu)$ of $n$ such that
$\lambda$ has distinct parts and the number of occurrences of the
part of $\mu$ in both $\lambda$ and $\mu$ is at least two.

\begin{kst}
\item[Case 1:] For a rooted partition $(\lambda,\,\mu)$ with
almost distinct parts but not with distinct parts, if $\lambda$
contains the part of $\mu$, then move this part from $\lambda$ to
$\mu$.

\item[Case 2:] For a rooted partition $(\lambda,\,\mu)$ with
almost distinct parts but not with distinct parts, if $\lambda$
does not contain the part of $\mu$, then move this part from $\mu$
to $\lambda$.
\end{kst}

It is easy to check that the above mapping is an involution.
Moreover, $\tau$ changes the parity of the root size. \qed

For example, there are nine  rooted partitions of 4 with almost
distinct parts:
$$\bar{4},\,\bar{3}+1,\,3+\bar{1},\,\bar{2}+2,\,\bar{2}+\bar{2},
\,2+\bar{1}+1,\,2+\bar{1}+\bar{1},\,1+\bar{1}+\bar{1}+\bar{1},\,\bar{1}+\bar{1}+\bar{1}+\bar{1}.$$
Applying the above involution, we get the following involution:
$$\bar{2}+2 \leftrightarrows
\bar{2}+\bar{2},\,2+\bar{1}+1\leftrightarrows 2+\bar{1}+\bar{1},\,
1+\bar{1}+\bar{1}+\bar{1}\leftrightarrows\bar{1}+\bar{1}+\bar{1}+\bar{1}.$$
The above involution does not apply to rooted partitions with
distinct parts: $\bar{4},\,\bar{3}+1,\,3+\bar{1}$.

The following correspondence can be regarded as a rooted partition
analogue of Euler's theorem.

\begin{thm} \label{o-2}
The number of the rooted partitions of $n$ into almost distinct
parts with odd root size   equals to the number of the rooted
partitions of $n$ with odd parts.
\end{thm}

\pf  We employ  Sylvester's bijection to construct a map from the
set of rooted partitions of $n$ into almost distinct parts with
odd root size to the set of rooted partitions of $n$ with odd
parts.

{\bf The map $\sigma$}: For a rooted partition $(\lambda,\,\mu)$
into almost distinct parts with odd root size, we apply the
inverse map of Sylvester's bijection $\varphi^{-1}$ to $\lambda$
to generate a partition $\alpha$ with odd parts.  Let $\beta$ be
the conjugate of $\mu$ which is a partition with equal odd parts.
Therefore $(\alpha,\, \beta)$ is a rooted partition with odd
parts.

{\bf The inverse map $\sigma^{-1}$}: For a rooted partition
$(\alpha,\,\beta)$ with odd parts, we apply   Sylvester's
bijection $\varphi$ to $\alpha$ to generate a  partition $\lambda$
 with distinct parts. Let $\mu$ be conjugate of $\beta$,
which is a partition into  equal parts with odd length. Thus
$(\lambda,\,\mu)$ is a rooted partition into almost distinct parts
with odd root size.

From Sylvester's bijection, one sees that $\sigma$ is also a
bijection. \qed

For example, there are  six rooted partitions of 4 into almost
distinct parts with odd root size:
$$\bar{4},\,\bar{3}+1,\,3+\bar{1},\,\bar{2}+2,\,2+\bar{1}+1,\,
1+\bar{1}+\bar{1}+\bar{1},$$ and there are six  rooted partitions
of 4 with odd parts:
$$\bar{3}+1,\,3+\bar{1},\,\bar{1}+1+1+1,\,
\bar{1}+\bar{1}+1+1,\,\bar{1}+\bar{1}+\bar{1}+1,
\bar{1}+\bar{1}+\bar{1}+\bar{1}.$$ Using that above bijection, we
have the following correspondence:
\begin{align*}
\begin{array}{ccc}
\bar{4} \leftrightarrows
\bar{1}+\bar{1}+\bar{1}+\bar{1}&\bar{3}+1\leftrightarrows
\bar{1}+\bar{1}+\bar{1}+1&
3+\bar{1}\leftrightarrows\bar{1}+1+1+1\\[6pt]
\bar{2}+2\leftrightarrows\bar{1}+\bar{1}+1+1&
2+\bar{1}+1\leftrightarrows
3+\bar{1}&1+\bar{1}+\bar{1}+\bar{1}\leftrightarrows\bar{3}+1.
\end{array}
\end{align*}

From the above Theorems \ref{o-1} and \ref{o-2}, we obtain Theorem
\ref{lem2} which serves as a combinatorial setting for Ramanujan's
identity \eqref{eqn2}. For Ramanujan's identity (\ref{eqn1}), we
need the following partition identity which also follows from the
above two theorems:

\begin{thm}\label{lem1}
The number of  rooted partitions  of $n$ with almost distinct
parts plus the number of  rooted partitions  of $n$ with distinct
parts is twice the number of  rooted partitions of $n$ with odd
parts.
\end{thm}

We now make a connection between rooted partitions with distinct
parts and odd parts and ordinary partitions with distinct parts
and odd parts. Chapman \cite{cha02} has shown that the series
(\ref{eqn4})
\[\sum_{d=1}^{\infty}q^d\prod_{n\neq
d}(1+q^{n}) \] is the generating function for ordinary partitions
with distinct parts with multiplicities (or weight) being their
lengths. Note that  the above series is also the generating
function for rooted partitions with distinct parts. This
generating function identity implies that there should be a
combinatorial correspondence between rooted partitions and
ordinary partitions with distinct parts.

 In fact, a simple correspondence goes as follows:
From a partition $\alpha$ with distinct parts, we can get
$l(\alpha)$ distinct rooted partitions $(\lambda,\,\mu)$ with
distinct parts by designating any  part of $\alpha$ as the part of
$\mu$ and keeping the remaining parts of $\alpha$ as  parts of
$\lambda$. This map is clearly reversible.

For example, there are two partitions of 4 with distinct parts:
$4,\,3+1$. The sum of their lengths is three, whereas there are
three rooted partitions of 4 with distinct parts:
$\bar{4},\,\bar{3}+1,\,3+\bar{1}.$

Thus we have the following theorem  on the relationship between
rooted partitions with distinct parts and  partitions with
distinct parts.

\begin{thm}\label{lem3}
The number of the rooted partitions of $n$ with distinct parts
equals  the sum of  lengths over  the partitions of $n$ with
distinct parts.
\end{thm}

Chapman \cite{cha02} has shown that identity \eqref{eqn5} is also
the generating function for the sum of the lengths of partitions
with odd parts:
\begin{align*}
\frac{1}{(q;q^2)_{\infty}}
\sum_{d=0}^{\infty}\frac{q^{2d+1}}{1-q^{2d+1}}
&=\sum_{d=0}^{\infty}\frac{1}{(q;q^2)_{d}(q^{2d+3};q^2)_{\infty}}\cdot\frac{q^{2d+1}}{(1-q^{2d+1})^2}\\[5pt]
&=\sum_{d=0}^{\infty}\sum_{m=1}^{\infty}\frac{mq^{(2d+1)m}}{(q;q^2)_{d}(q^{2d+3};q^2)_{\infty}}\\[5pt]
&=\sum_{d=0}^{\infty}\sum_{\lambda\in
O }n_{\lambda }(2d+1)q^{|\lambda|}\\[5pt]
&=\sum_{\lambda\in O }l(\lambda)q^{|\lambda|}.
\end{align*}

Using the formulation of rooted partitions with odd parts and the
above generating function, we obtain the following relation
between rooted partitions and ordinary partitions, and we give a
combinatorial proof of this fact. Theorem \ref{lem3} and the
following Theorem \ref{lem4} will be necessary to transform the
formulations of Ramanujan's identities with rooted partitions to
combinatorial settings with ordinary partitions.

\begin{thm}\label{lem4}
The number of rooted partitions of $n$ with odd parts  equals the
sum of  lengths over the partitions  of $n$ with odd parts.
\end{thm}

\pf  In fact, for a partition $\beta$ of $n$ with odd parts, we
may  get $l(\beta)$ distinct  rooted partitions $(\lambda,\,\mu)$
of $n$ with odd parts by designating any part of $\beta$ as the
parts of $\mu$ and keep the remaining parts of $\beta$ as parts of
$\lambda$. Assume that $d$ is a part that appears $m$ times
$(m\geq 2)$ in $\beta$. Then we may choose $\mu$ as the partition
with $d$ repeated $i$ times, where $i=1, 2, \ldots, m$. \qed

For example, there are two partitions of 4  with odd parts namely
$3+1,\,1+1+1+1$, the sum of lengths is six. For rooted partitions
of 4,
 we see
that there are also six rooted partitions with odd parts:
$\bar{3}+1,\,3+\bar{1},\,\bar{1}+1+1+1,\,
\bar{1}+\bar{1}+1+1,\,\bar{1}+\bar{1}+\bar{1}+1,
\bar{1}+\bar{1}+\bar{1}+\bar{1}.$

\section{ Ramanujan's Identities }

In this section, we will reduce  Ramanujan's identities
(\ref{eqn1}) and (\ref{eqn2}) to the two weighted forms
(\ref{we1}) and (\ref{we2}) of Euler's theorem.  The left hand
sides of (\ref{eqn1}) and  (\ref{eqn2}) have partition
interpretations as given by Andrews \cite{and86}. The first
summations on the right hand sides of \eqref{eqn1} and
\eqref{eqn2} can be interpreted combinatorially in term of
ordinary partitions with multiplicities as given by Theorem
\ref{lem2} and \ref{lem1}. The second summations on the right hand
sides of \eqref{eqn1} and \eqref{eqn2} have partition
interpretations in terms of the rank.

 Combining
Theorems \ref{lem3} and \ref{lem4} on the relations between rooted
partitions and ordinary partitions, we may transform Theorem
\ref{lem1} on rooted partitions to following statement on ordinary
partitions:

\begin{thm} \label{t4.1}
The number of rooted partitions  of $n$ with almost distinct parts
equals  the sum of twice  the lengths over partitions of $n$ with
odd parts minus  the sum of lengths over partitions of $n$ with
distinct parts. In terms of generating functions, we have
$$(-q;q)_{\infty}\sum_{d=1}^{\infty}
\frac{q^d}{1-q^d}=\sum_{\lambda \in
O}2l(\lambda)q^{|\lambda|}-\sum_{\mu \in D}l(\mu)q^{|\mu|}.$$
\end{thm}

We proceed to demonstrate that with the aid of the above theorem,
Ramanujan's identity \eqref{eqn1} can be restated as the weighted
form (\ref{we1}) of Euler's theorem.  Using the following relation
due to Andrews \cite{and86}:
\begin{equation*}
\sum_{n=0}^{\infty}\left[ (-q;q)_{\infty}-(-q;q)_n
\right]=\sum_{n=1}^{\infty}nq^n(1+q)(1+q^2)\cdots(1+q^{n-1}),
\end{equation*}  the left side of
Ramanujan's identity \eqref{eqn1} equals the generating function
for the sum of the largest parts over the partitions with distinct
parts:
$$\sum_{n=0}^{\infty}\left[
(-q;q)_{\infty}-(-q;q)_n \right]=\sum_{\mu \in D}
\mu_1~q^{|\mu|}.$$ It is easy to see that the second summation on
the right hand of (\ref{eqn1}), that is,
$$1+\sum_{n=1}^{\infty}\frac{q^{n+1\choose 2}}{(-q;q)_n},$$
equals the generating function for partitions into distinct parts
with even rank minus the generating function for  partitions into
distinct parts with odd rank. Note that the coefficient of $q^m$
in the above series has been studied in \cite{and88}. Therefore,
we have
 \begin{equation}
 \label{g-z}
 -\frac{1}{2}(-q;q)_{\infty}+\frac{1}{2}
\left[1+\sum_{n=1}^{\infty}\frac{q^{n+1\choose
2}}{(-q;q)_n}\right] =-\sum_{\stackrel{\mu \in D}{ r(\mu)\
\mbox{\scriptsize odd}}}q^{|\mu|}.
\end{equation}

From the above interpretations and Theorem \ref{t4.1}, one sees
the right side of Ramanujan's identity \eqref{eqn1} equals the
generating function for the sum of twice the lengths over
partitions with odd parts minus the generating function for the
sum of  lengths over partitions with distinct parts minus the
generating function for partitions into distinct parts with odd
rank:
\begin{align*}
& (-q;q)_{\infty} \left[ -\frac{1}{2}+\sum_{d=1}^{\infty}
\frac{q^d}{1-q^d} \right]+\frac{1}{2}\left[1+\sum_{n=1}^{\infty}
\frac{q^{n+1\choose 2}}{(-q;q)_n}\right]\\[8pt]
&=\sum_{\lambda \in O}2l(\lambda)q^{|\lambda|}-\sum_{\mu \in
D}l(\mu)q^{|\mu|}-\sum_{\stackrel{\mu \in D}{ r(\mu)\
\mbox{\scriptsize odd}}}q^{|\mu|}.
\end{align*}
We now reach the conclusion that Ramanujan's identity \eqref{eqn1}
can be restated as the weighted form \eqref{we1} of Euler's
theorem:
\begin{equation}\label{eqn9} \sum_{\mu \in
D}\left(\mu_1+l(\mu)+\frac{1-(-1)^{r(\mu)}}{2}\right)q^{|\mu|}
=\sum_{\lambda \in O}2l(\lambda)q^{|\lambda|}.
 \end{equation}
Thus, we have obtained a combinatorial proof of (\ref{eqn1}) based
on a weighted form of Euler's theorem.

Similarly,   combining Theorems \ref{lem3}  and \ref{lem4} on the
relations between rooted partitions and ordinary partitions, we
may transform Theorem \ref{lem2} on rooted partitions to the
following assertion for ordinary partitions:

\begin{thm}\label{t4.2}
The number of  rooted partitions  of $n$ into almost distinct
parts with even length equals  the sum of  lengths over the
partitions of $n$ into odd parts  minus the sum of  lengths over
partitions of $n$ into distinct parts. In terms of generating
functions, we have
$$(-q;q)_{\infty}\sum_{d=1}^{\infty}\frac{q^{2d}}{1-q^{2d}}
=\sum_{\lambda \in O}l(\lambda)q^{|\lambda|}-\sum_{\mu \in
D}l(\mu)q^{|\mu|}.$$
\end{thm}

Using the following relation due to Andrews \cite{and86}:
\begin{align*}
&\hskip-2mm\sum_{n=0}^{\infty}\left[\frac{1}{(q;q^2)_{\infty}}
-\frac{1}{(q;q^2)_n}\right]\\[5pt]
&=\sum_{n=0}^{\infty}\frac{nq^{2n+1}}{(1-q)(1-q^3)\cdots(1-q^{2n+1})},
\end{align*}
the left side of Ramanujan's identity \eqref{eqn2} equals the
generating function of  the sum of half of its largest part minus
one over the partitions into odd parts:
$$\sum_{n=0}^{\infty}\left[\frac{1}{(q;q^2)_{\infty}}
-\frac{1}{(q;q^2)_n}\right]=\sum_{\lambda \in
O}\frac{\lambda_1-1}{2}~q^{|\lambda|}.$$
By using the above
relation (\ref{g-z}) and Theorem \ref{t4.2}, one sees that the
right side of Ramanujan's identity \eqref{eqn2} equals the
generating function for the sum of  lengths over partitions into
odd parts  minus the generating function for the sum of lengths
over partitions into distinct parts minus the generating function
for partitions into distinct parts with odd rank:
\begin{align*}
&(-q;q)_{\infty}\left[-\frac{1}{2}+
\sum_{d=1}^{\infty}\frac{q^{2d}}{1-q^{2d}}\right]
+\frac{1}{2}\left[1+\sum_{n=1}^{\infty}\frac{q^{n+1\choose 2}}{(-q;q)_n}\right]\\[5pt]
&=\sum_{\lambda \in O}l(\lambda)q^{|\lambda|}-\sum_{\mu \in
D}l(\mu)q^{|\mu|}-\sum_{\stackrel{\mu \in D}{ r(\mu)\
\mbox{\scriptsize odd}}}q^{|\mu|}.
\end{align*}

So Ramanujan's identity \eqref{eqn2} can be recast as the weighted
form \eqref{we2} of Euler's theorem:
\begin{equation}\label{eqn6}
\sum_{\mu \in
D}\left(l(\mu)+\frac{1-(-1)^{r(\mu)}}{2}\right)q^{|\mu|}=\sum_{\lambda
\in O}\left(l(\lambda)-\frac{\lambda_1-1}{2}\right)q^{|\lambda|}.
\end{equation}

\vspace{6mm}

\noindent{\bf Acknowledgments.} We would like to thank George E.
Andrews for valuable comments. This work was supported by the 973
Project on Mathematical Mechanization, the National Science
Foundation, the Ministry of Education, and the Ministry of Science
and Technology of China.

\end{document}